\documentclass{article}
\usepackage{amsmath,amsfonts,amssymb,array,tabularx,enumerate}
\usepackage{epigraph}
\usepackage{enumitem}      
\usepackage{arxiv}
\usepackage{float}
\usepackage[utf8]{inputenc} 
\usepackage{fvextra} 
\usepackage{capt-of}   
\usepackage[T1]{fontenc}    
\usepackage{verbatim}
\DeclareUnicodeCharacter{03B1}{\ensuremath{\alpha}}
\usepackage{hyperref}       
\usepackage{url}  
\usepackage{booktabs}
\usepackage{amsfonts}  
\usepackage{nicefrac}
\usepackage{microtype}  
\usepackage{lipsum} 
\usepackage{graphicx}
\usepackage{doi}
\usepackage[numbers,sort&compress]{natbib}
\usepackage{tikz}
\usepackage{ifthen}
\usetikzlibrary{arrows.meta,positioning}
\usepackage[RGB]{xcolor} 
\definecolor{gold}{RGB}{255,215,0} 
\newtheorem{thm}{Theorem}[section]
\newtheorem{cor}[thm]{Corollary}
\newtheorem{lem}[thm]{Lemma}
\newtheorem{prop}[thm]{Proposition}

\newtheorem{rema}[thm]{Remark}
\newtheorem{defi}[thm]{Definition}

\title{A Classification of Elements of Sequence Space $Seq(\mathbb{R})$}
\author{Mohsen Soltanifar\thanks{Contact: mohsen.soltanifar[at]alumni.utoronto.ca.} \href{https://orcid.org/0000-0002-5989-0082}{\includegraphics[scale=0.06]{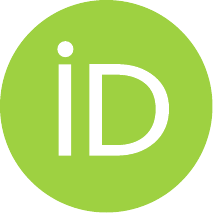}}
\\ Biostatistics Division, Dalla Lana School of Public Health, University of Toronto\\
620-155 College Street, Toronto, ON M5T 3M7, Canada}
 
\renewcommand{\shorttitle}{\textit{arXiv} Template}

\begin{document}
\maketitle

\begin{abstract}
The sequence space of all real-valued sequences, denoted $Seq(\mathbb{R})$, is typically investigated through the lens of infinite-dimensional vector spaces, utilizing Banach space norms or Schauder bases. This work proposes a complementary, constructive classification based instead on the asymptotic limit profile encoded by the pair $(\liminf a_n, \limsup a_n)$. We demonstrate that this perspective naturally partitions $Seq(\mathbb{R})$ into seven mutually disjoint macroscale blocks, covering behaviors from finite convergence to bounded and unbounded oscillation. For each block, we provide explicit closed-form representative sequences and establish that every constituent class possesses the cardinality of the continuum. Furthermore, we investigate the structural relationships between these blocks at two distinct levels of granularity. At the macroscale, we employ injective mappings to define an idealized connectivity graph, while at the microscale, we introduce a connection relation governed by the Hadamard (pointwise) product. This dual analysis reveals a rich directed graph structure where the block of finite convergent sequences functions as both the only subspace and as a global attractor with no outgoing connections. Statistical comparisons between the idealized and realized adjacency matrices indicate that the pointwise product structure realizes approximately two-thirds of the theoretically possible macroscale relations. Ultimately, this partition-based framework endows the seemingly chaotic space $Seq(\mathbb{R})$ with a transparent, geometrically interpretable internal structure.
\end{abstract}

\vspace{1em}

\epigraph{``Mathematics is the science of the infinite, its goal the symbolic comprehension of the infinite with human, that is finite, means.''}{--- Hermann K.H. Weyl (1885–-1955)}

\section{Introduction}\label{sec1.}
\subsection{Real Valued Sequences}\label{sec1.1}
The systematic study of real-valued sequences emerged as a cornerstone of modern real analysis during the nineteenth century. While infinite processes had appeared implicitly in the work of Newton and Leibniz, it was the shift toward rigor, initiated by Cauchy and later clarified by Weierstrass and others, that placed sequences of real numbers at the heart of the theory of limits and convergence (see, e.g., \cite{Grabiner1981,RogersBoman2013}). Sequences provided a flexible language to express approximation, to formulate the Cauchy criterion, and to capture completeness properties of the real line. Historical studies of this transition document how the step from intuitive infinitesimals to $\varepsilon$--$\delta$ arguments was largely mediated by sequences and series \cite{Grabiner1981,RogersBoman2013}.\par 

In contemporary real analysis, real-valued sequences and their series are treated as one of the primary objects of study, alongside real numbers and real-valued functions. Standard texts typically devote early chapters to basic notions such as boundedness, monotonicity, subsequences, and Cauchy sequences, and then develop fundamental theorems including the Bolzano--Weierstrass theorem and the equivalence between Cauchy and convergent sequences in $\mathbb{R}$ \cite{Tao2016,Hunter2014,Deshpande2004,LokuBraha2024}. These results not only underpin the rigorous development of calculus but also guide later topics such as series, metric spaces, and functional analysis. Introductory treatments by authors such as Tao, Hunter, Deshpande, Loku and Braha, and others all emphasize sequences as the natural setting in which students first encounter the interplay between algebraic structure, order, and completeness \cite{Tao2016,Hunter2014,Deshpande2004,LokuBraha2024}.\par 

\subsection{Motivation}\label{sec1.2}
The sequence space of all real-valued sequences on the natural numbers, often denoted $Seq(\mathbb{R})$ or $\mathbb{R}^{\mathbb{N}}$, is an infinite-dimensional vector space with an extremely rich internal structure. As a set of functions $a\colon\mathbb{N}\to\mathbb{R}$, it contains the familiar convergent and Cauchy sequences that encode limits and completeness, but it also hosts a vast collection of divergent, oscillatory, and highly irregular sequences. From the standpoint of cardinality, $Seq(\mathbb{R})$ has the size of the continuum and admits Hamel bases that are necessarily uncountable and nonconstructive; from the analytic standpoint, familiar Banach sequence spaces such as $\ell^p$, $\ell^\infty$, and $c_0$ arise as distinguished linear subspaces equipped with norms and (often) countable Schauder bases \cite{Deshpande2004,Tao2016}. For the well-known Banach subspaces of the sequence space $Seq(\mathbb{R})$ we have the nested proper inclusion chain ($1 \leq p<+\infty$): $c_{00} \subset \ell^{1} \subset \ell^{2} \subset \cdots \subset \ell^{p} \subset c_{0}
\subset c \subset \ell^{\infty} \subset Seq(\mathbb{R}).$\par 

Because of this breadth, many classical approaches focus on particular ``regular'' subspaces or impose additional structural conditions: one may restrict attention to bounded or monotone sequences, Cauchy sequences, or sequences belonging to $\ell^p$ or $c_0$; alternatively, one may classify sequences through summability methods (such as Ces\`aro or Abel summation), through topological properties in metric or product topologies, or through measure-theoretic or probabilistic considerations in the study of random sequences \cite{Deshpande2004,Hunter2014,LokuBraha2024}. These viewpoints have been extremely successful, but they typically either narrow the focus to specific well-behaved classes or rely on infinite expansions with respect to a chosen basis.\par 

The present work proposes a complementary, constructive description of $Seq(\mathbb{R})$ that does not depend on a particular Hamel or Schauder basis. This approach has been previously applied in the case of function space $F(\mathbb{R},\mathbb{R})$ \cite{Soltanifar2023}. Instead, we organize sequences according to their asymptotic ``limit profile,'' encoded by the pair $(\liminf a_n,\limsup a_n)$, and show that this perspective leads to a finite partition of $Seq(\mathbb{R})$ into seven macroscale blocks. Within and between these blocks we then investigate finer ``connection'' relations that reflect how sequences can be transformed into one another while preserving or modifying their limit behaviour. In this way, the enormous and seemingly chaotic space $Seq(\mathbb{R})$ acquires a more transparent, geometrically interpretable structure that complements existing classification schemes based on topology, summability, or basis representations.\par 

\subsection{Study Outline}\label{sec1.3}
This paper is organized into three main sections. In Section~\ref{sec2} we collect the requisite background from set theory, linear algebra, and the theory of special classes of sequences, providing a common framework for the subsequent analysis. Section~\ref{sec3} develops the central results: we construct a finite partition of the sequence space \(Seq(\mathbb{R})\) into seven blocks, supply explicit representatives for each block, and then analyze the relationships between these blocks at both a macroscale and a microscale level using a directed graph viewpoint. Finally, Section~\ref{sec4} offers a discussion of the main findings, highlights the structural insights gained from this partition-based perspective, and outlines several directions for future work.\par 

\section{Preliminaries}\label{sec2}
The reader who has studied linear algebra and real analysis is well equipped with the following set of definition, propositions, theorems and remarks regarding infinite sequences.\par 

\subsection{Set-Theoretic Foundations}
\label{sec2.1}

In this subsection we recall basic set-theoretic notions and standard properties of real-valued sequences that will be used throughout the paper. Our goal is to fix notation and to highlight the central role of the limit inferior and limit superior in describing the long–run behaviour of a sequence.\par 

\begin{defi}[Sequence space \(Seq(\mathbb{R})\)]
\label{def2.1}
We denote by
\begin{eqnarray}
Seq(\mathbb{R}) &:=& \{ a = (a_n)_{n \geq 1} : a_n \in \mathbb{R} \text{ for all } n \in \mathbb{N} \}
\end{eqnarray}
the set of all real-valued sequences indexed by the natural numbers.  When convenient we write \(a \in Seq(\mathbb{R})\) as \(a = (a_n)\) or \(a = (a_n)_{n \geq 1}\).
\end{defi}

\begin{defi}[Limit profile of the sequence]
\label{def2.2}
Let \(a = (a_n) \in  Seq(\mathbb{R})\). For each \(n \in \mathbb{N}\) define the tail infimum and tail supremum
\begin{eqnarray}
  \alpha_n &:=& \inf_{k \geq n} a_k, \\
  \beta_n  &:=& \sup_{k \geq n} a_k.
\end{eqnarray}
The \emph{limit inferior} and \emph{limit superior} of \(a\) are given by
\begin{eqnarray}
  \liminf_{n \to \infty} a_n &:=& \sup_{n \in \mathbb{N}} \alpha_n, \label{eq:def-liminf}\\
  \limsup_{n \to \infty} a_n &:=& \inf_{n \in \mathbb{N}} \beta_n, \label{eq:def-limsup}
\end{eqnarray}
with values in the extended real line \(\overline{\mathbb{R}} := \mathbb{R} \cup \{-\infty,+\infty\}\).
We will frequently write
\begin{eqnarray}
  L_1(a) &:=& \liminf_{n \to \infty} a_n, \\
  L_2(a) &:=& \limsup_{n \to \infty} a_n.
\end{eqnarray}
We refer to the pair $(L_1(a),L_2(a))$ as the limit profile of the sequence $a.$
\end{defi}

\begin{prop}[General properties of \(\liminf\) and \(\limsup\)]
\label{prop2.3}
Let \(a = (a_n) \in Seq(\mathbb{R})\) and define \(\alpha_n,\beta_n\) as in Definition~\ref{def2.2}. Then:
\begin{enumerate}[label=(\roman*)]
  \item The sequence \((\alpha_n)_{n \in \mathbb{N}}\) is non-decreasing and the sequence \((\beta_n)_{n \in \mathbb{N}}\) is non-increasing.
  \item For every \(n \in \mathbb{N}\) we have
  \begin{eqnarray}
    \alpha_n &\leq& a_n \;\leq\; \beta_n.
  \end{eqnarray}
  \item The limits
  \begin{eqnarray}
    \lim_{n \to \infty} \alpha_n &=& \liminf_{n \to \infty} a_n, \\
    \lim_{n \to \infty} \beta_n  &=& \limsup_{n \to \infty} a_n
  \end{eqnarray}
  always exist in \(\overline{\mathbb{R}}\).
  \item We always have the inequality
  \begin{eqnarray}
    \liminf_{n \to \infty} a_n &\leq& \limsup_{n \to \infty} a_n.
  \end{eqnarray}
\end{enumerate}
\end{prop}

\begin{cor}\label{col2.4}
Let $a = (a_n)$ be a real-valued sequence and let $\alpha$ be a real number. Define the shifted sequence by  $(a + \alpha)_n = a_n + \alpha$. Then, for both
the liminf and limsup we have
$L_i(a + \alpha) = L_i(a) + \alpha, \qquad i = 1,2.$
\end{cor} 

\begin{thm}[Relationship between \(\liminf\), \(\limsup\), and the usual limit]
\label{thm2.4}
Let \(a = (a_n) \in Seq(\mathbb{R})\) and set
\begin{eqnarray}
  L_1(a) &:=& \liminf_{n \to \infty} a_n, \\
  L_2(a) &:=& \limsup_{n \to \infty} a_n.
\end{eqnarray}
Then the following statements hold:
\begin{enumerate}[label=(\roman*)]
  \item If the (finite) limit \(\displaystyle \lim_{n \to \infty} a_n\) exists and is equal to \(L \in \mathbb{R}\), then
  \begin{eqnarray}
    L_1(a) &=& L_2(a) \;=\; L(a).
  \end{eqnarray}
  \item Conversely, if \(L_1(a) = L_2(a) = L(a) \in \mathbb{R}\), then \((a_n)\) converges to \(L(a)\), that is,
  \begin{eqnarray}
    \lim_{n \to \infty} a_n &=& L(a).
  \end{eqnarray}
  \item If \(L_1(a) < L_2(a)\), then the sequence \((a_n)\) does not converge in \(\mathbb{R}\); in this case \((a_n)\) is either divergent to \(\pm\infty\) or oscillatory between at least two distinct cluster values.
\end{enumerate}
\end{thm}

\begin{rema}
\label{rema2.5}
In the extended real setting, the values \(L_1(a)\) and \(L_2(a)\) may be equal to \(\pm\infty\). For example, the sequence \(a_n = n\) satisfies $\liminf_{n \to \infty} a_n=
\limsup_{n \to \infty} a_n \;=\; +\infty,$ while \(a_n = (-1)^n n\) has $ \liminf_{n \to \infty} a_n = -\infty, $ and $  \limsup_{n \to \infty} a_n =  +\infty.$ These extremal behaviours will play a role in our later classification of \(Seq(\mathbb{R})\) into seven macroscale blocks.
\end{rema}

\subsection{Linear-Algebraic Foundations}\label{sec2.2}

We now recall the linear structure of the sequence space \(Seq(\mathbb{R})\) and discuss two complementary notions of basis: the algebraic (Hamel) basis and the topological (Schauder) basis. The former captures the purely set-theoretic size of \(Seq(\mathbb{R})\), while the latter reflects more analytic information in classical Banach sequence spaces.

\begin{defi}[Vector space structure on \(Seq(\mathbb{R})\)]
\label{def2.6}
We regard \(Seq(\mathbb{R})\) as a real vector space under pointwise operations: for \(a = (a_n)\), \(b = (b_n)\) in \(Seq(\mathbb{R})\) and \(\lambda \in \mathbb{R}\) we define
\begin{eqnarray}
  (a + b)_n &:=& a_n + b_n, \\
  (\lambda a)_n &:=& \lambda a_n,
\end{eqnarray}
for all \(n \in \mathbb{N}\).
With these operations, \(Seq(\mathbb{R})\) is an infinite-dimensional real vector space.
\end{defi}

\begin{defi}[Hamel basis]
\label{def2.7}
Let \(X\) be a real vector space. A subset \(B \subseteq X\) is called a \emph{Hamel basis} (or algebraic basis) of \(X\) if:
\begin{enumerate}[label=(\roman*)]
  \item The elements of \(B\) are linearly independent.
  \item Every element \(x \in X\) can be written as a \emph{finite} linear combination of elements of \(B\); that is,
  \begin{eqnarray}
    x &=& \sum_{k=1}^{m} \lambda_k b_k,
  \end{eqnarray}
  for some \(m \in \mathbb{N}\), scalars \(\lambda_1,\dots,\lambda_m \in \mathbb{R}\), and distinct \(b_1,\dots,b_m \in B\).
\end{enumerate}
If such a set \(B\) exists, the cardinality \(|B|\) is called the (Hamel) dimension of \(X\).
\end{defi}

\begin{rema}[Cardinality of a Hamel basis of \(Seq(\mathbb{R})\)]
\label{rem2.8}
It is a classical result in set-theoretic linear algebra (assuming the Axiom of Choice) that the Hamel dimension of \(Seq(\mathbb{R})\) as a vector space over \(\mathbb{R}\) is equal to the cardinality of the continuum:
\begin{eqnarray}
\dim_{H}\bigl( Seq(\mathbb{R}) \bigr)
  &=& \mathfrak{c}.
\end{eqnarray}
In particular, no countable Hamel basis of \(Seq(\mathbb{R})\) exists. From the point of view of concrete analysis, any Hamel basis for \(Seq(\mathbb{R})\) is therefore necessarily highly non-constructive.
\end{rema}

\begin{defi}[Schauder basis]
\label{def2.9}
Let \(X\) be a Banach space over \(\mathbb{R}\). A sequence \((e_n)_{n \geq 1} \subseteq X\) is called a \emph{Schauder basis} of \(X\) if for every \(x \in X\) there exists a unique sequence of scalars \((\alpha_n)_{n \geq 1} \subseteq \mathbb{R}\) such that
\begin{eqnarray}
  x &=& \lim_{m \to \infty} \sum_{n=1}^{m} \alpha_n e_n,
\end{eqnarray}
with convergence taken in the norm of \(X\). The coefficients \((\alpha_n)\) are called the (Schauder) coordinates of \(x\) with respect to the basis \((e_n)\).
\end{defi}

\begin{rema}[Schauder bases in classical sequence spaces]
\label{rema2.10}
In standard Banach sequence spaces such as \(\ell^p\) (\(1 \leq p < \infty\)) and \(c_0\), the canonical unit vectors
\begin{eqnarray}
  e^{(k)} &:=& (0,0,\dots,0,\underbrace{1}_{k\text{th position}},0,\dots),
\end{eqnarray}
for \(k \in \mathbb{N}\), form a countable Schauder basis. Thus every element of these spaces can be represented as a convergent infinite linear combination of the \(e^{(k)}\).

By contrast, the full sequence space \(Seq(\mathbb{R})\) equipped with its natural product topology is not a Banach space and does not admit such a simple Schauder basis. For our purposes, the key point is the contrast between:
\begin{enumerate}[label=(\roman*)]
  \item the \emph{algebraic} viewpoint, where \(Seq(\mathbb{R})\) has a Hamel basis of cardinality \(\mathfrak{c}\), and
  \item the \emph{analytic} viewpoint, where familiar sequence spaces admit countable Schauder bases describing their elements through convergent series.
\end{enumerate}
Our classification of \(Seq(\mathbb{R})\) into seven macroscale blocks will be constructive and finitary in spirit, and will be independent of any particular choice of Hamel or Schauder basis.
\end{rema}

\subsection{Special Sequences}\label{sec2.3}

\begin{lem}\label{lem2.1}
Let $a=(a_n)_{n\ge 1}$ be a real sequence with infinitely many zeros. Then $L_1(a)$ and $L_2(a)$ can never both be $+\infty$ nor both be $-\infty$. In particular, $a$ can
never diverge to $\pm\infty$.
\end{lem}
\textbf{Proof.}
This is straightforward consequence of the fact that $-\infty \leq L_1(a) \leq 0 \leq L_2(a) \leq +\infty.$ \\
$\Box$

\begin{lem}\label{lem2.2}
Let $a=(a_n)_{n\ge 1}$ be a non–zero real sequence which is not
eventually zero (i.e.\ $a_n\neq 0$ for infinitely many indices $n$).
Then we can construct connectors $c=(c_n)_{n\ge 1}$ such that, for the
Hadamard (pointwise) product  ~\cite{diDioLanger2025}:
\[
(a\odot c)_n := a_n c_n \qquad (n\ge 1),
\]
the following patterns occur:
\begin{itemize}
  \item[(A)] $a\odot c$ has $L_1(a\odot c) = -\infty$ and $L_2(a\odot c)$ finite;
  \item[(B)] $a\odot c$ has finite $L_1(a\odot c) < L_2(a\odot c)$;
  \item[(C)] $a\odot c$ has finite $L_1(a\odot c)$ and $L_2(a\odot c) = +\infty$;
  \item[(D)] $a\odot c$ has $L_1(a\odot c) = -\infty$ and $L_2(a\odot c) = +\infty$.
\end{itemize}
Here $L_1$ and $L_2$ denote, respectively, the $\liminf$ and $\limsup$
of the corresponding sequence.
\end{lem}
\textbf{Proof.}
Let $a=(a_n)_{n\ge1}$ be as in the statement and put
\[
S:=\{n\in\mathbb{N}:a_n\neq 0\}.
\]
By hypothesis $S$ is infinite; fix an enumeration $S=\{s_1,s_2,\dots\}$.
For any connector $c=(c_n)_{n\ge1}$ we write $b=a\odot c$ for the
Hadamard product, i.e.\ $b_n=(a\odot c)_n=a_n c_n$.

\smallskip\noindent
\emph{(A–pattern).}
Define $c$ by
\[
c_{s_k}:=-\frac{k}{a_{s_k}},\qquad n\notin S\ \Longrightarrow\ c_n:=0 .
\]
Then $b_{s_k}=-k$ for all $k$ and $b_n=0$ for $n\notin S$.
Hence $b$ has a subsequence equal to $-k\to-\infty$ and infinitely many
zeros, so
\[
L_1(b)=\liminf_{n\to\infty}b_n=-\infty,
\qquad
L_2(b)=\limsup_{n\to\infty}b_n=0.
\]

\smallskip\noindent
\emph{(B–pattern).}
Split $S$ into two infinite subsequences, for instance
\[
S^+ := \{s_{2k}:k\ge1\},\qquad
S^- := \{s_{2k+1}:k\ge0\}.
\]
Define $b$ by
\[
b_n:=(a\odot c)_n :=
\begin{cases}
  1,  & n\in S^+,\\[2pt]
 -1,  & n\in S^-,\\[2pt]
  0,  & n\notin S,
\end{cases}
\]
and then set $c_n:=b_n/a_n$ whenever $a_n\neq0$, and $c_n:=0$ when
$a_n=0$.  By construction $b$ takes the values $\pm1$ infinitely often
(and $0$ possibly as well), hence
\[
L_1(b)=-1,\qquad L_2(b)=1,
\]
so $L_1(b)$ and $L_2(b)$ are finite with $L_1(b)<L_2(b)$.

\smallskip\noindent
\emph{(C–pattern).}
Choose any infinite subsequence $\{s_k\}\subset S$ and define
\[
b_{s_k}:=k>0\quad (k\ge1),\qquad b_n:=0 \ \text{for } n\notin\{s_k\}_{k\ge1},
\]
and again put $c_n:=b_n/a_n$ when $a_n\neq0$, $c_n:=0$ when $a_n=0$.
Then $b_{s_k}=k\to+\infty$ while $b_n=0$ for infinitely many $n$, so
\[
L_1(b)=0,\qquad L_2(b)=+\infty.
\]

\smallskip\noindent
\emph{(D–pattern).}
Using the same partition $S=S^+\cup S^-$ as above, define
\[
b_{s_{2k}}:=k,\qquad b_{s_{2k+1}}:=-k\quad (k\ge1),
\qquad b_n:=0\ \text{for } n\notin S,
\]
and set $c_n:=b_n/a_n$ for $a_n\neq0$, $c_n:=0$ when $a_n=0$.
Then $b$ has subsequences $k\to+\infty$ and $-k\to-\infty$, hence
\[
L_1(b)=-\infty,\qquad L_2(b)=+\infty.
\]

In each case we have explicitly constructed a connector $c$ such that
the Hadamard product $a\odot c$ exhibits the desired pair
$(L_1,L_2)$, which completes the proof.\\
$\Box$

\section{Main Results}\label{sec3}
\subsection{Partition of $Seq(\mathbb{R})$ with Scenario Classification \& Examples}\label{sec3.1}

In this section we adopt a constructive viewpoint on the sequence space \(Seq(\mathbb{R})\) by organizing it into three successive themes: (i) the existence of certain canonical “blocks,” (ii) the explicit construction of these blocks, and (iii) the assessment of their size. We begin by partitioning \(Seq(\mathbb{R})\) into a finite collection of blocks obtained from the interaction of the two classical notions \(\liminf\) and \(\limsup\): two real-valued infinite sequences lie in the same block if and only if they share the same pair \((\liminf, \limsup)\), as made precise below. We then turn, in Theorem~\ref{thm3.2}, to the question of representation, and ask whether each such block admits at least one explicit infinite sequence in closed form that can serve as its representative. Finally, still within the framework of Theorem~\ref{thm3.2}, we address the cardinality of each block and determine the size of these classes inside the ambient space \(Seq(\mathbb{R})\).\par 

\begin{defi}[Asymptotic Equivalence]\label{Def3.1}
Let $a,b\in Seq(\mathbb{R}),$ with associated limits $L_1(a),L_2(a)$ and $L_1(b),L_2(b),$ respectively. Then, $a$ is equivalent to  $b$  asymptotically, denoted by  $a R_{asymp} b$, whenever (i) $L_1(a)=L_1(b),$ (ii) $L_2(a)=L_2(b).$
\end{defi}

As is easily verified, the relation \(\mathcal{R}_{\mathrm{asymp}}\) introduced in Definition~\ref{Def3.1} is an equivalence relation on \(Seq(\mathbb{R})\), and therefore it induces a natural partition of the sequence space into its corresponding equivalence classes.\par 

\begin{thm}[Finite Partition by Limit Profile]\label{thm3.2}
The sequence space \(Seq(\mathbb{R})\) can be partitioned into seven pairwise disjoint blocks, each of cardinality equal to the continuum \(\mathfrak{c}\).
\end{thm}
\textbf{Proof.} 
Let \(a \in Seq(\mathbb{R})\) with \(L_1(a) \leq L_2(a)\) and \(L_2(a) \in [-\infty,+\infty]\). We distinguish three general situations according to the value of \(L_2(a)\). First, if \(L_2(a)=+\infty\), then \(L_1(a)\) may equal \(+\infty\), be a finite real number, or equal \(-\infty\) (3 possibilities). Second, if \(-\infty < L_2(a) < +\infty\), then \(L_1(a)\) may coincide with \(L_2(a)\), be a finite real strictly smaller than \(L_2(a)\), or equal \(-\infty\) (again 3 possibilities). Third, if \(L_2(a)=-\infty\), then necessarily \(L_1(a)=-\infty\) (1 possibility). Altogether this yields \(3+3+1=7\) possible configurations, and hence 7 corresponding blocks. These blocks, together with one explicit representative sequence for each, are listed in Table~\ref{tab1}. Moreover, a straightforward perturbation of each representative sequence (e.g., adding a constant $\alpha$ in $(0,1)$ as in Corollary~\ref{col2.4}) shows that the cardinality of every block is the continuum \(\mathfrak{c}\).\par
\begin{figure}[H] 
\centering 
\includegraphics[clip,width=0.75\columnwidth]{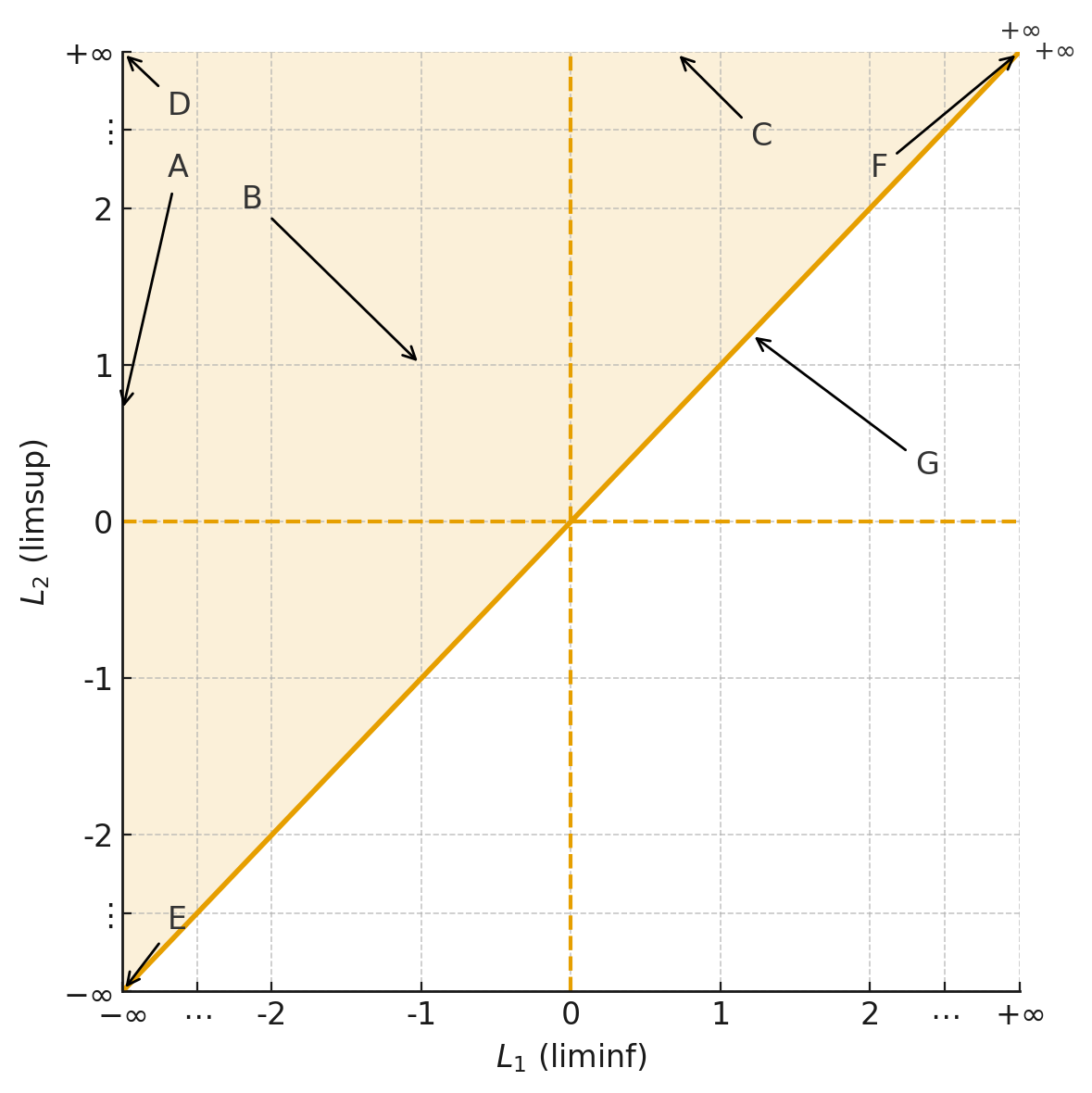}
\caption{Partition of $Seq(\mathbb{R})$ by the $(L_1,L_2)$-plane, where $L_1(a)=\liminf a_n$ and $L_2(a)=\limsup a_n$, into seven admissible regions labelled $A,\dots,G$ satisfying $L_1 \le L_2$. Each region corresponds to a distinct asymptotic behaviour class of real-valued sequences (finite convergent, finite–finite divergent, and the various one-sided or two-sided infinite cases).\label{fig1}}
\end{figure}

 
\begin{table}[htbp]\label{tab1}
\centering
\caption{List of seven representatives blocks of partition of $\operatorname{Seq}(\mathbb{R})$ with associated representative sequence and size of the block.}
\begin{tabular}{cclllc}
\hline
\# & Block & Definition of $(L_{1}(a),L_{2}(a))$ &Description & Representative $(a_n)$ & Size \\ \hline
1 & $A$ & $L_{1}(a)=-\infty,\ -\infty<L_{2}(a)<+\infty$ 
  &lower-unbounded oscillatory & $ n\bigl(\sin(\tfrac{n\pi}{2})-1\bigr),$ 
  & continuum $(\mathfrak{c})$ \\[2pt]
2 & $B$ & $-\infty<L_{1}(a)<L_{2}(a)<+\infty$ 
  &bounded oscillatory & $\sin(\tfrac{n\pi}{2}),$ 
  & continuum $(\mathfrak{c})$ \\[2pt]
3 & $C$ & $-\infty<L_{1}(a)<L_{2}(a)=+\infty$ 
  &upper-unbounded oscillatory & $ n\bigl(\sin(\tfrac{n\pi}{2})+1\bigr),$ 
  & continuum $(\mathfrak{c})$ \\[2pt]
4 & $D$ & $L_{1}(a)=-\infty,\ L_{2}(a)=+\infty$ 
  &fully-unbounded oscillatory & $n\sin(\tfrac{n\pi}{2}),$ 
  & continuum $(\mathfrak{c})$ \\[2pt]
5 & $E$ & $L_{1}(a)=-\infty,\ L_{2}(a)=-\infty$ 
  &downward divergent & $-n,$ 
  & continuum $(\mathfrak{c})$ \\[2pt]
6 & $F$ & $L_{1}(a)=+\infty,\ L_{2}(a)=+\infty$ 
  &upward divergent & $n,$ 
  & continuum $(\mathfrak{c})$ \\[2pt]
7 & $G$ & $-\infty<L_{1}(a)=L_{2}(a)<+\infty$ 
  &finite convergent & $\frac{1}{n},$ 
  & continuum $(\mathfrak{c}).$ \\ \hline
\end{tabular}
\end{table}
$\Box$\\
Figure~\ref{fig1} presents the seven blocks introduced in the Theorem\ref{thm3.2}. Here, all spaces
$c_{00},\ell^{1},\ell^{2},\cdots, \ell^{p},c_{0}$ ( $1\leq p<+\infty$)
correspond to the origin point $(L_{1},L_{2}) = (0,0)$, since every sequence in these spaces
converges to $0$. The space $c$ of all convergent sequences occupies the diagonal line $L_{1} = L_{2}$, representing convergence to a finite limit $L \in \mathbb{R}$. The space $\ell^{\infty}$ of bounded sequences fills the \emph{finite region} of the plane where both $L_{1}$ and $L_{2}$ are finite and satisfy $L_{1} \leq L_{2}$. Finally, the ambient space $Seq(\mathbb{R})$ covers the entire admissible area of Figure~\ref{fig1}, including all divergent and unbounded cases.

The following corollary is immediate from Theorem~\ref{thm3.2}:
\begin{cor}\label{cor3.3}
Let $Seq(\mathbb{R})$ be the vector space of all real sequences with given partition as in Theorem~\ref{thm3.2}. Then: 
\begin{enumerate}[label=(\roman*)]
  \item $G$ is a linear subspace of $Seq(\mathbb{R})$.
  \item None of the blocks $A,B,C,D,E,F$ is a linear subspace of $Seq(\mathbb{R})$.
  \item The Hamel dimension of $G$ equals the cardinality of the continuum, \(\dim_{H} G = \mathfrak{c}\).
\end{enumerate}
\end{cor}
\textbf{Proof.} (i) Let $\lambda\in\mathbb{R}$ and $a,b\in G.$ Then, $|L(\lambda.a)|=|\lambda||L(a)|<+\infty$ and $|L(a+b)|\leq |L(a)|+|L(b)|<+\infty,$ implying $\lambda.a\in G$ and $a+b\in G,$ respectively. (ii) Let $X\neq G,$ and $a\in X.$ Then, $0.a=0\notin X.$ Hence, $X$ is not closed under scalar product. (iii) Given $G = c_{0} \oplus \operatorname{span}\{(1)_{n\geq1}\}$ it follows that $\dim_H G = \dim_H c_0 + 1 = \mathfrak{c} + 1 = \mathfrak{c}.$\\
$\Box$

We conclude this section by showing that $\ell^\infty=B \cup G$ is the  maximal proper subspace of $Seq(\mathbb{R})$ with respect to inclusion among such block--subspaces.

\begin{lem}\label{lem3.4}
Let $Seq(\mathbb{R})$ be partitioned into the seven blocks $A,\dots,G$,
and let $M \subset Seq(\mathbb{R})$ be a union of some of these blocks.
If $X$ is one of the blocks with $M \cap X \neq \emptyset$, then
$X \subset M$.
\end{lem}

\textbf{proof.}
By assumption $M = \bigcup_{i\in I} B_i$ for some
$\{B_i : i\in I\} \subset \{A,B,C,D,E,F,G\}$, where the blocks form a
partition of $Seq(\mathbb{R})$.  Let $x \in M \cap X$.  Then $x \in M$,
so $x \in B_j$ for some $j \in I$, and $x \in X$ by choice.
Since the blocks are pairwise disjoint, $X \cap B_j \neq \emptyset$
implies $X = B_j$, hence $X \subset M$.\par 
$\Box$

\begin{thm}[Block subspaces and inclusion-based maximality]\label{thm3.5}
Let $Seq(\mathbb{R})$ be partitioned into the seven blocks $A,\dots,G$
as in Theorem~\ref{thm3.2}, and let $M \subset Seq(\mathbb{R})$
be a union of some of these blocks.  Then $M$ is a linear subspace of
$Seq(\mathbb{R})$ if and only if
\[
  M \in \bigl\{\,G,\; B \cup G,\; Seq(\mathbb{R})\,\bigr\}.
\]
\end{thm}
\textbf{Proof.}
Since the seven blocks form a partition of $Seq(\mathbb{R})$, every
union of blocks is of the form
\[
  M = \bigcup_{i \in I} B_i,
\]
where $\{B_i : i \in I\} \subset \{A,B,C,D,E,F,G\}$.  Any linear
subspace of $Seq(\mathbb{R})$ must contain the zero sequence, and the
zero sequence belongs to the block $G$ (its liminf and limsup are both
$0$).  Hence every block--union that is a subspace necessarily contains
$G$, so we may write
\[
M = G \cup X, \qquad X \subset \{A,B,C,D,E,F\}.
\]
We complete the proof in the following three steps:\par 
\medskip\noindent
\textbf{Step 1: Symmetry under sign change.}
The map $a \mapsto -a$ sends
\[
  A \longleftrightarrow C, \qquad
  E \longleftrightarrow F, \qquad
  B \longrightarrow B, \qquad
  D \longrightarrow D, \qquad
  G \longrightarrow G,
\]
because $(L_1(-a),L_2(-a)) = (-L_2(a),-L_1(a))$.  Since any subspace is
closed under multiplication by $-1$, the presence of $A$ forces the
presence of $C$, and the presence of $E$ forces the presence of $F$.
Thus, in a block--union subspace, $A$ and $C$ either occur together or
not at all, and similarly for $E$ and $F$.\par 
\medskip\noindent
\textbf{Step 2: Adding unbounded blocks forces all blocks.}
We now show that if $X$ contains at least one of the unbounded blocks
$A,C,D,E,F$, then necessarily $M = Seq(\mathbb{R})$.  This is done by frequent application of Lemma~\ref{lem3.4} in the following four cases.\par 
\emph{Case 1: $D \subset M$.}
First, consider the two sequences
\[
  d^{(1)}_n := (-1)^n n^2, \qquad
  d^{(2)}_n := (-1)^{n+1} n^2 - n.
\]
Both $d^{(1)}$ and $d^{(2)}$ have subsequences diverging to $+\infty$
and to $-\infty$, hence $d^{(1)},d^{(2)} \in D \subset M$.  Their sum
is $d^{(1)}_n + d^{(2)}_n = -n$, which satisfies
$L_1 = L_2 = -\infty$ and therefore lies in block $E$.  Because $M$ is
a union of blocks, this implies $E \subset M$, and by Step~1 we also
obtain $F \subset M$.

Next, using $E,F\subset M$ take the sequences
\[
  e_n := -n \in E, \qquad
  f_n := n + \sin\frac{n\pi}{2} \in F.
\]
Then $e+f = (\sin(\tfrac{n\pi}{2}))_{n\ge 1}$ is bounded and oscillatory
with limit inferior $-1$ and limit superior $1$, so $e+f \in B$.  Since
$E,F \subset M$ and $M$ is a union of blocks, this gives $B \subset M$.

Finally, consider
\[
  d_n := (-1)^n n \in D, \qquad
  e_n := -n \in E, \qquad
  f_n := n \in F.
\]
Define $a := d+e$ and $c := d+f$.  Then
\[
  a_n =
  \begin{cases}
    -2n, & n \text{ odd},\\[2pt]
    0,   & n \text{ even},
  \end{cases}
  \qquad
  c_n =
  \begin{cases}
    0,   & n \text{ odd},\\[2pt]
    2n,  & n \text{ even}.
  \end{cases}
\]
The sequence $a$ is unbounded below and bounded above by $0$, so
$L_1(a) = -\infty$ and $L_2(a) = 0$, hence $a \in A$.  The sequence
$c$ is unbounded above and bounded below by $0$, so $L_1(c) = 0$ and
$L_2(c) = +\infty$, hence $c \in C$.  Since $d,e,f \in M$, both $a$
and $c$ belong to $M$, and the union-of-blocks assumption implies
$A,C \subset M$.

We have therefore shown that $D \subset M$ implies
\[
  \{A,B,C,D,E,F\} \subset M,
\]
so $M = Seq(\mathbb{R})$.

\smallskip\noindent
\emph{Case 2: $A \subset M$ or $C \subset M$.}
By Step~1, $A \subset M$ implies $C \subset M$, and conversely.  Let
$a^{(A)}$ and $c^{(C)}$ be the representative sequences of blocks $A$
and $C$, respectively, given by:
\[
  a^{(A)}_n = n\!\left(\sin\frac{n\pi}{2}-1\right), \qquad
  c^{(C)}_n = n\!\left(\sin\frac{n\pi}{2}+1\right).
\]
Then, their sum given by
\[
  d_n := a^{(A)}_n + c^{(C)}_n = 2n \sin\frac{n\pi}{2}
\]
satisfies $d_{4k+1} = 2(4k+1) \to +\infty$ and
$d_{4k+3} = -2(4k+3) \to -\infty$ as $k \to \infty$.  Hence
$L_1(d) = -\infty$ and $L_2(d) = +\infty$, so $d \in D$.  Since $M$ is
a union of blocks and $d \in M$, we obtain $D \subset M$.  Case~1
above will then imply $M = Seq(\mathbb{R})$.

\smallskip\noindent
\emph{Case 3: $E \subset M$ or $F \subset M$.}
Suppose first that $E \subset M$. By Step 1, $E\subset M$ implies $F\subset M$, and conversely. Thus $E,F \subset M$(If $F \subset M$, the same argument implies $E \subset M$.).

Now define two sequences $e = (e_n)$ and $f = (f_n)$ by
\[
  e_n :=
  \begin{cases}
    -n^2, & n \text{ odd},\\[2pt]
    -n,   & n \text{ even},
  \end{cases}
  \qquad
  f_n :=
  \begin{cases}
     n,   & n \text{ odd},\\[2pt]
     n^2, & n \text{ even}.
  \end{cases}
\]
We have $e_n \to -\infty$ and $f_n \to +\infty$ as $n \to \infty$, so
$e \in E$ and $f \in F$.  Their sum $d := e+f$ satisfies
$d_{2k+1} = - (2k+1)^2 + (2k+1) \to -\infty$ on odd indices and
$d_{2k} = -2k + (2k)^2 \to +\infty$ on even indices.  Hence
$L_1(d) = -\infty$ and $L_2(d) = +\infty$, so $d \in D$.  Because
$E,F \subset M$ and $M$ is a union of blocks, this yields $D \subset M$. Case~1
above will then imply $M = Seq(\mathbb{R})$.\par 

\emph{Case 4: $B \subset M$.}
If $B \subset M$ and $M$ also contains at least one of the unbounded
blocks $A,C,D,E,F$, then Cases~1, 2, and~3 show that
$M = Seq(\mathbb{R})$.  Thus, in order for $M$ to be a proper
subspace, the presence of $B$ must be the only possible extension
beyond $G$, and we shall analyse this case separately in Step~3.

Combining Cases~1, 2, and~4, and the remark in Case~3, we conclude that
\[
  X \cap \{A,C,D,E,F\} \neq \varnothing
  \quad\Longrightarrow\quad
  M = Seq(\mathbb{R}).
\]

\medskip\noindent
\textbf{Step 3: Bounded cases and maximality.}
It remains to analyse the situation where $M$ is a proper subset of
$Seq(\mathbb{R})$.  By Step~2, this is only possible if
\[
  X \subset \{B\},
\]
so either $X = \varnothing$ or $X = \{B\}$.

If $X = \varnothing$, then $M = G$, which is the space $c$ of
convergent sequences.  By Corollary~\ref{cor3.3}, $G$ is a
linear subspace of $Seq(\mathbb{R})$.

If $X = \{B\}$, then
\[
  M = B \cup G.
\]
Thus $M$ is exactly the space $\ell^\infty$ of bounded sequences which is a classical linear subspace of $Seq(\mathbb{R})$.\par 
$\Box$

\begin{rema}[Schauder bases for the three block spaces]
In the Banach-space sense, among the three canonical block spaces only
$M = G = c$ admits a Schauder basis: for example, the family
$\{f_0,f_1,f_2,\dots\}$ with $f_0 = (1,1,1,\dots)$ and $f_n = e_n$ for
$n \ge 1$ is a Schauder basis of $c$ in $\|\cdot\|_\infty$.
The bounded space $M = B \cup G = \ell^\infty$ is nonseparable and
therefore has no Schauder basis.
The full sequence space $M = Seq(\mathbb{R})$ is not a Banach space
under the supremum norm (only $\ell^\infty$ is), so a Schauder-basis
question is not canonically posed there without choosing a different
topology.
\end{rema}

\subsection{The Relationship between the  Blocks}\label{sec3.2}

In this section we examine how the seven blocks \(A,B,C,D,E,F,G\) of the sequence space \(Seq(\mathbb{R})\) are related at two complementary levels. At a macroscale level, we exploit the fact that all blocks have the same cardinality to construct explicit injections between them, mapping elements of one block into another and vice versa. At a microscale level, we introduce the notion of \emph{connection}, which we use, whenever possible, to represent elements in a given block as pointwise products of elements drawn from other blocks.\par

\subsubsection{Macroscale}\label{sec3.2.1}

\begin{lem}[Injectivity of the coding map]
Let $S:= Seq(\mathbb{R})$ be the set of all real-valued sequences.
Define a bijection $\sigma:\mathbb{R}\to(0,1)$ by
\begin{eqnarray}
\sigma(x) &:=& \frac12\left(1+\frac{x}{1+|x|}\right), \qquad x\in\mathbb{R}.
\end{eqnarray}
Let $w_n:=2^{-2^n}$ for $n\ge1$. These weights satisfy
$w_n>\sum_{k>n}w_k$ for every $n$.
For $a=(a_n)_{n\ge1}\in S$ set
\begin{eqnarray}
\Phi(a) &:=& \sum_{n=1}^{+\infty} w_n \sigma(a_n)
          \;=\; \sum_{n=1}^{+\infty} 2^{-2^n}\,\sigma(a_n)\in(0,1),
\end{eqnarray}
and denote
\begin{eqnarray}
c(a) &:=& \Phi(a).
\end{eqnarray}
Then the map $\Phi:S\to(0,1)$ is injective. In particular, the real number
$c(a)\in(0,1)$ is a unique code for the sequence $a$.
\end{lem}

\textbf{Proof.}
Suppose $\Phi(a)=\Phi(b)$, that is,
\begin{eqnarray}
\sum_{n=1}^{+\infty} w_n \sigma(a_n) &=& \sum_{n=1}^{+\infty} w_n \sigma(b_n).
\end{eqnarray}
If for some $m$ we had $\sigma(a_m)\neq\sigma(b_m)$, choose the least such
$m$. Then the difference of the two series at index $m$ has magnitude at
least $w_m$, while the tail $\sum_{k>m}w_k$ is strictly smaller than $w_m$;
this is impossible. Hence $\sigma(a_n)=\sigma(b_n)$ for all $n\geq 1$, so
$a_n=b_n$ for all $n\geq 1$ because $\sigma$ is bijective. Thus $\Phi$ is
injective, and the real number $c(a)\in(0,1)$ is a unique code for the
sequence $a$.\\
$\Box$
\begin{thm}[Macroscale connection between the seven blocks]\label{thm3.4}
Let $S = Seq(\mathbb{R})$ be the set of all real-valued sequences, and for  $a= (a_n)_{n\ge1}\in S$ set $L_1(a),L_2(a)$ as defined above with seven blocks $(A,B,C,D,E,F,G)$ as in Table~\ref{tab1}. Then for any two distinct blocks $X,Y \in \{A,B,C,D,E,F,G\}$ there exists an injective map:
\begin{eqnarray}\label{inj-map}
T_{(Y,X)} : Y \longrightarrow X.
\end{eqnarray}
\end{thm}
\textbf{Proof.} We discuss the existence of injective maps $T_X's$ in three steps as follows: \par 

\noindent\textbf{Step 1: Injective maps $T_X$ from the complement to $X$.}
For each target region $X\in\{A,B,C,D,E,F,G\}$ we define a map:
\begin{eqnarray}\label{com-inj-map}
T_X &:& S\setminus X \longrightarrow X  \\
T_X((a_n)_{n\geq1}) &=& (T_X(a)_n)_{n\geq1}\nonumber
\end{eqnarray}
using only the code $c(a)$ and the index $n$. We list the definitions and
their liminf/limsup:

\begin{itemize}
\item[(i)] Target $G$ ($-\infty<L_1=L_2<+\infty$).
  \[
  T_G(a)_n := c(a)\quad\text{for all }n.
  \]
  Then $T_G(a)$ is constant, hence in $G$.

\item[(ii)] Target $F$ ($L_1=L_2=+\infty$).
  \[
  T_F(a)_n := n + c(a).
  \]
  Clearly $T_F(a)_n\to+\infty$, so $T_F(a)\in F$.

\item[(iii)] Target $E$ ($L_1=L_2=-\infty$).
  \[
  T_E(a)_n := -n - c(a),
  \]
  so $T_E(a)_n\to-\infty$ and $T_E(a)\in E$.

\item[(iv)] Target $D$ ($L_1=-\infty,\ L_2=+\infty$).
  \[
  T_D(a)_{2n} := n + c(a),\qquad
  T_D(a)_{2n-1} := -n - c(a).
  \]
  Along even indices the sequence tends to $+\infty$; along odd indices it
  tends to $-\infty$, so $T_D(a)\in D$.

\item[(v)] Target $C$ ($-\infty<L_1<L_2=+\infty$).
  \[
  T_C(a)_{2n}   := n,\qquad
  T_C(a)_{2n-1} := c(a)+\frac1n.
  \]
  Then $T_C(a)_{2n-1}\to c(a)$ and $T_C(a)_{2n}\to+\infty$, so
  $L_1(T_C(a))=c(a)$ and $L_2(T_C(a))=+\infty$, hence $T_C(a)\in C$.

\item[(vi)] Target $B$ ( $-\infty<L_1<L_2<+\infty$).
  Define
  \[
  \alpha(a):=c(a),\qquad \beta(a):=c(a)+1,
  \]
  and set
  \[
  T_B(a)_{2n-1} := \alpha(a)+\frac1n,\qquad
  T_B(a)_{2n}   := \beta(a)-\frac1n.
  \]
  Then $L_1(T_B(a))=\alpha(a)$ and $L_2(T_B(a))=\beta(a)$ with
  $\alpha(a)<\beta(a)$, so $T_B(a)\in B$.

\item[(vii)] Target $A$ ($L_1=-\infty,\ -\infty<L_2<+\infty$).
  \[
  T_A(a)_{2n}   := -n,\qquad
  T_A(a)_{2n-1} := c(a)-\frac1n.
  \]
  Then $T_A(a)_{2n}\to-\infty$ and $T_A(a)_{2n-1}\to c(a)$, so
  $L_1(T_A(a))=-\infty$ and $L_2(T_A(a))=c(a)$; hence $T_A(a)\in A$.
\end{itemize}

\smallskip
\noindent\textbf{Step 2: Each $T_X$ is injective.}
In each case, the image sequence determines $c(a)$ uniquely:

\begin{itemize}
\item For $T_G$, $c(a)=T_G(a)_n$ for any $n$.
\item For $T_F$, $c(a)=T_F(a)_n-n$.
\item For $T_E$, $c(a)=-T_E(a)_n-n$.
\item For $T_D$, $c(a)=\lim_{n\to+\infty}\bigl(T_D(a)_{2n}-n\bigr)$.
\item For $T_C$, $c(a)=\lim_{n\to+\infty}T_C(a)_{2n-1}$.
\item For $T_B$, $c(a)=\liminf_{n\to+\infty}T_B(a)_n$.
\item For $T_A$, $c(a)=\limsup_{n\to+\infty}T_A(a)_n$.
\end{itemize}
Thus if $T_X(a)=T_X(b)$, then the corresponding codes coincide:
$c(a)=c(b)$. Since $\Phi$ is injective, $c(a)=c(b)$ implies $a=b$. Hence
each $T_X$ is a one-to-one map.

\noindent\textbf{Step 3: Injective maps $T_{(Y,X)}$ from  $Y$ to $X$.} Given arguments in Step 1 and Step 2, it now sufficient to consider:
\begin{eqnarray}\label{inj-map2}
T_{Y,X} &:& Y \longrightarrow X.\\
T_{Y,X} (a) &=& (T_{X}|Y) (a).\nonumber 
\end{eqnarray}

Finally, summarizing all implications, the adjacency matrix
$U=[u_{XY}]$ (rows = source $X$, columns = target $Y$) for distinctive pairs $(X,Y)$ is
\[
u_{XY} =
\begin{cases}
1, & \text{if } X \to Y,\\
0, & \text{otherwise},
\end{cases}
\qquad X,Y\in\{A,B,C,D,E,F,G\},
\]
and concretely:
\begin{eqnarray}\label{matrix-u}
U=
\begin{array}{c|rrrrrrr}
   & A & B & C & D & E & F & G\\ \hline
A & 0 & 1 & 1 & 1 & 1 & 1 & 1\\
B & 1 & 0 & 1 & 1 & 1 & 1 & 1\\
C & 1 & 1 & 0 & 1 & 1 & 1 & 1\\
D & 1 & 1 & 1 & 0 & 1 & 1 & 1\\
E & 1 & 1 & 1 & 1 & 0 & 1 & 1\\
F & 1 & 1 & 1 & 1 & 1 & 0 & 1\\
G & 1 & 1 & 1 & 1 & 1 & 1 & 0
\end{array}.    
\end{eqnarray}

$\Box$
\begin{rema}\label{rema3.5}
The total number of constructive injections between seven blocks in Theorem~\ref{thm3.4} is 42 where each of them is a restriction of one of seven  constructive injections in definition (\ref{com-inj-map}) to one of six blocks $Y$ distinctive from the target block $X.$ Figure~\ref{fig2} presents sample case for $Y=F$ and $X=G.$
\begin{figure}[H] 
\centering 
\includegraphics[clip,width=0.75\columnwidth]{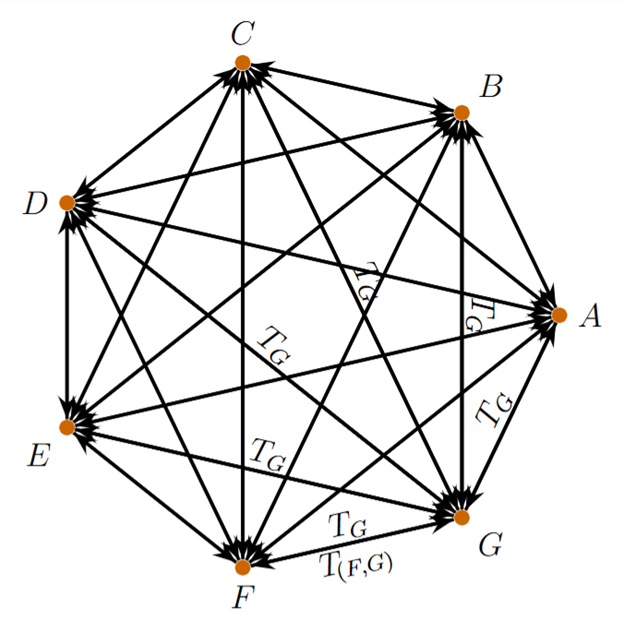}
\caption{Heptagon with with macroscale directed connectivity graph between seven blocks (vertices) $A,B,C,D,E,F,G$ and directed edges from $A,\dots,F$ to $G$. Each edge is tagged by the global map $T_G$, while the edge $F\to G$ also represents the specific map $T_{(F,G)}$.\label{fig2}}
\end{figure}

\end{rema}
\begin{rema}\label{rema3.6}
Given any two distinct blocks $X,Y \in \{A,B,C,D,E,F,G\}$. Then, by two applications of Theorem~\ref{thm3.4} we have: $|Y|=|T_{(Y,X)}(Y)|\leq |X|=|T_{(X,Y)}(X)|\leq |Y|.$ Consequently:
\begin{eqnarray}\label{card-block}
|A|=|B|=|C|=|D|=|E|=|F|=|G|.    
\end{eqnarray}
\end{rema}
\subsubsection{Microscale}\label{sec3.2.2}

\begin{defi}[Connection of Sequences]\label{Def3.7}
Let \(a,b \in Seq(\mathbb{R})\). We say that \(a\) is \emph{connected} to \(b\), and write
\(a \stackrel{\mathrm{conn}}{\sim}  b\), if there exists a sequence \(c \in Seq(\mathbb{R})\) such that
\(a \odot c \in [b]\), where the pointwise (Hadamard) product is defined by:
\begin{eqnarray}
(a \odot c)_n = a_n c_n,\qquad n \ge 1.    
\end{eqnarray}
\end{defi}

\begin{rema}\label{rema3.8}
If in Definition~\ref{Def3.7} we take \(c\) to be the constant sequence
\(1 = (1,1,\ldots)\), then \(a \odot c = a\) and the relation
\(\stackrel{\mathrm{conn}}{\sim} \) reduces to the equivalence relation
\(R_{\mathrm{asymp}}\) introduced in Definition~\ref{Def3.1}. In this sense,
\(\stackrel{\mathrm{conn}}{\sim} \) extends the asymptotic equivalence relation \(R_{\mathrm{asymp}}\).
\end{rema}

\begin{rema}\label{rema3.9}
Definition~\ref{Def3.7} naturally induces a notion of connection between the blocks at the microscale. For any \(X,Y \in \{A,B,C,D,E,F,G\}\), we say that \(X\) is \emph{connected} to \(Y\), and write \(X \to Y\), if for every sequence \(a \in X\) there exists a sequence \(c \in Seq(\mathbb{R})\) such that the pointwise Hadamard product \(a \odot c\) belongs to \(Y\).
\end{rema}

\begin{thm}[Microscale connection between the seven blocks]\label{thm3.10}
Let $S = Seq(\mathbb{R})$ be the set of all real-valued sequences with  seven blocks $(A,B,C,D,E,F,G)$ as in Table~\ref{tab1}. Then any two distinct blocks $X,Y \in \{A,B,C,D,E,F,G\}$ are connected at microscale level.
\end{thm}
\textbf{Proof.} 
\noindent\textbf{Step 1: Case Discussion} 
Recall that for $X,Y\in\{A,B,C,D,E,F,G\}$ we write $X\to Y$ if for every $a\in X$ there exists a connector $c\in Seq(\mathbb{R})$ such that the Hadamard product $a\odot c$ lies in $Y$.  Our goal is to determine, for each ordered pair of distinct blocks $(X,Y)$, whether $X\to Y$ holds.\par 
We shall proceed by \emph{fixing the target block $Y$} and then examining all six possible sources $X\neq Y$.\par 

\medskip
\noindent\textbf{Case $Y=G$ (finite limit).}
By definition
\[
G=\bigl\{a\in Seq(\mathbb{R}): -\infty<L_1(a)=L_2(a)<+\infty\bigr\}.
\]
Given any $a\in Seq(\mathbb{R})$, choose the constant connector $c_n\equiv 0$.
Then $(a\odot c)_n=0$ for all $n$, so $a\odot c$ is the constant zero
sequence and hence
\[
L_1(a\odot c)=L_2(a\odot c)=0,
\]
that is, $a\odot c\in G$.  Since $a$ was arbitrary, we obtain
$X\to G$ for every $X\in\{A,B,C,D,E,F,G\}$.  Restricting to distinct
blocks, all six implications $X\to G$ with $X\neq G$ are true.

\medskip
\noindent\textbf{Case $Y=F$ ($L_1=L_2=+\infty$).}
Here
\[
F=\bigl\{a\in Seq(\mathbb{R}): L_1(a)=L_2(a)=+\infty\bigr\};
\]
we need the product $a\odot c$ to diverge to $+\infty$.

\smallskip\noindent
\emph{Positive result.}
If $a\in E$, then $a_n\to -\infty$, so $a_n<0$ eventually.  With the
constant connector $c_n\equiv -1$ we have $a\odot c=-a$, and therefore
\[
L_1(a\odot c)=L_2(a\odot c)=+\infty,
\]
so $a\odot c\in F$ for every $a\in E$.  Thus $E\to F$.

\smallskip\noindent
\emph{Negative results.}
By Lemma~\ref{lem2.1}, if a sequence has infinitely many zeros then no Hadamard product with it can belong to $F$ (its liminf and limsup cannot both be $+\infty$).  Hence, if a block $X$ contains even one sequence with infinitely many zeros, then $X\nrightarrow F$: taking that particular $a\in X$ witnesses the failure of the definition of $X\to F$.\par 
We can choose such a sequence in each of the following blocks:
\begin{itemize}
  \item in $A$: $a_n=-n$ for odd $n$, $a_n=0$ for even $n$.  Then
        $L_1(a)=-\infty$, $L_2(a)=0$, and $a$ has infinitely many zeros;
  \item in $B$: $a_n=0$ for even $n$, $a_n=1$ for odd $n$, so
        $L_1(a)=0$, $L_2(a)=1$;
  \item in $C$: $a_n=n$ for odd $n$, $a_n=0$ for even $n$, so
        $L_1(a)=0$, $L_2(a)=+\infty$;
  \item in $D$: for example $a_{3k-2}=-k$, $a_{3k-1}=0$, $a_{3k}=k$
        ($k\ge1$), so $L_1(a)=-\infty$, $L_2(a)=+\infty$;
  \item in $G$: the constant zero sequence $a_n=0$ has $L_1(a)=L_2(a)=0$.
\end{itemize}
In each case $a$ has infinitely many zeros, so by
Lemma~\ref{lem2.1} no connector $c$ can produce
$L_1(a\odot c)=L_2(a\odot c)=+\infty$.  Thus
\[
A\nrightarrow F,\quad B\nrightarrow F,\quad
C\nrightarrow F,\quad D\nrightarrow F,\quad G\nrightarrow F.
\]

Therefore, with target $Y=F$, the only true implication among the six distinct possibilities $X\to F$ is
\[
E\to F.
\]

\medskip
\noindent\textbf{Case $Y=E$ ($L_1=L_2=-\infty$).}
This is completely symmetric to the previous case.  We have
\[
E=\bigl\{a\in Seq(\mathbb{R}): L_1(a)=L_2(a)=-\infty\bigr\}.
\]

\smallskip\noindent
\emph{Positive result.}
If $a\in F$ (so $a_n\to +\infty$), the connector $c_n\equiv -1$ gives
$a\odot c=-a$, hence
\[
L_1(a\odot c)=L_2(a\odot c)=-\infty,
\]
and so $a\odot c\in E$ for every $a\in F$.  Thus $F\to E$.

\smallskip\noindent
\emph{Negative results.}
We reuse exactly the same “infinitely many zeros’’ witnesses in $A,B,C,D,G$
listed above.  For each such $a$ no Hadamard product $a\odot c$ can have
$L_1=L_2=-\infty$, again by Lemma~\ref{lem2.1}.
Consequently,
\[
A\nrightarrow E,\quad B\nrightarrow E,\quad
C\nrightarrow E,\quad D\nrightarrow E,\quad G\nrightarrow E.
\]

Thus, for target $Y=E$, the only true implication among the six distinct
cases $X\to E$ is
\[
F\to E.
\]

\medskip
\noindent\textbf{Case $Y=D$ ($L_1=-\infty$, $L_2=+\infty$).}
Here
\[
D=\bigl\{a\in Seq(\mathbb{R}): L_1(a)=-\infty,\ L_2(a)=+\infty\bigr\}.
\]
By the D–pattern in Lemma~\ref{lem2.2}, for every
$a\in X$ with $X\in\{A,B,C,D,E,F\}$ we can construct a connector $c$
such that $a\odot c$ satisfies
\[
L_1(a\odot c)=-\infty,\qquad L_2(a\odot c)=+\infty,
\]
hence $a\odot c\in D$.  Therefore,
\[
X\to D\qquad\text{for all }X\in\{A,B,C,D,E,F\}.
\]

For $X=G$ consider the zero sequence $a_n=0\in G$.  For any connector
$c$ the product $a\odot c$ is again identically zero, with
$L_1=L_2=0$, so it never lands in $D$.  Hence $G\nrightarrow D$.

Thus, with $Y=D$, five of the six implications $X\to D$ (those with
$X\neq G$) are true, and only $G\to D$ is false.

\medskip
\noindent\textbf{Case $Y=C$ ($-\infty<L_1<L_2=+\infty$).}
Now
\[
C=\bigl\{a\in Seq(\mathbb{R}): -\infty<L_1(a)<L_2(a)=+\infty\bigr\}.
\]
Using the C–pattern in Lemma~\ref{lem2.2}, we can, for
every $a\in X$ with $X\in\{A,B,C,D,E,F\}$, construct a connector $c$
such that
\[
L_1(a\odot c)=0,\qquad L_2(a\odot c)=+\infty,
\]
so that $a\odot c\in C$.  Hence
\[
X\to C\qquad\text{for all }X\in\{A,B,C,D,E,F\}.
\]

For $X=G$ we again take the zero sequence $a_n=0\in G$.  For any
connector $c$ the product $a\odot c$ is identically zero with finite
upper and lower limits, so $L_2(a\odot c)\neq +\infty$ and thus
$a\odot c\notin C$.  Therefore $G\nrightarrow C$.

So with $Y=C$ the only failure among the six implications $X\to C$ is
$G\to C$.

\medskip
\noindent\textbf{Case $Y=B$ (finite $L_1<L_2<+\infty$).}
We have
\[
B=\bigl\{a\in Seq(\mathbb{R}): -\infty<L_1(a)<L_2(a)<+\infty\bigr\}.
\]
By the B–pattern of Lemma~\ref{lem2.2}, for each
$a\in X$ with $X\in\{A,B,C,D,E,F\}$ we can construct a connector $c$
so that, for instance,
\[
L_1(a\odot c)=-1<L_2(a\odot c)=1,
\]
and hence $a\odot c\in B$.  Thus
\[
X\to B\qquad\text{for all }X\in\{A,B,C,D,E,F\}.
\]

For $X=G$ we again choose $a_n=0\in G$.  For any connector $c$,
$a\odot c$ has $L_1=L_2=0$, so $a\odot c\in G$ and never in $B$.
Therefore $G\nrightarrow B$.

With target $Y=B$ we therefore get five true implications $X\to B$
($X\neq G$) and one false implication $G\to B$.

\medskip
\noindent\textbf{Case $Y=A$ ($L_1=-\infty$, $L_2$ finite).}
Finally,
\[
A=\bigl\{a\in Seq(\mathbb{R}): L_1(a)=-\infty,\ L_2(a)<+\infty\bigr\}.
\]
Using the A–pattern in Lemma~\ref{lem2.2}, we can, for
every $a\in X$ with $X\in\{A,B,C,D,E,F\}$, build a connector $c$ such
that
\[
L_1(a\odot c)=-\infty,\qquad L_2(a\odot c)=0,
\]
so that $a\odot c\in A$.  Therefore
\[
X\to A\qquad\text{for all }X\in\{A,B,C,D,E,F\}.
\]
Restricting to distinct blocks, this yields
$X\to A$ for all $X\in\{B,C,D,E,F\}$.

For $X=G$, the zero sequence $a_n=0\in G$ once more blocks the
implication: for any connector $c$, $a\odot c\equiv 0$ has
$L_1=L_2=0$, so it never belongs to $A$.  Hence $G\nrightarrow A$.

Thus, with target $Y=A$, again five of the six implications $X\to A$
(with $X\neq A$) are true and $G\to A$ is false.\par 
\noindent\textbf{Step 2: Summary of Cases}
Summarizing all implications, the adjacency matrix
$V=[v_{XY}]$ (rows = source $X$, columns = target $Y$) for distinctive pairs $(X,Y)$ is
\[
v_{XY} =
\begin{cases}
1, & \text{if } X \to Y,\\
0, & \text{otherwise},
\end{cases}
\qquad X,Y\in\{A,B,C,D,E,F,G\},
\]
and concretely:
\begin{eqnarray}\label{matrix-v}
V=
\begin{array}{c|rrrrrrr}
   & A & B & C & D & E & F & G\\ \hline
A & 0 & 1 & 1 & 1 & 0 & 0 & 1\\
B & 1 & 0 & 1 & 1 & 0 & 0 & 1\\
C & 1 & 1 & 0 & 1 & 0 & 0 & 1\\
D & 1 & 1 & 1 & 0 & 0 & 0 & 1\\
E & 1 & 1 & 1 & 1 & 0 & 1 & 1\\
F & 1 & 1 & 1 & 1 & 1 & 0 & 1\\
G & 0 & 0 & 0 & 0 & 0 & 0 & 0
\end{array}.    
\end{eqnarray}
$\Box$

\begin{figure}[H] 
\centering 
\includegraphics[clip,width=0.75\columnwidth]{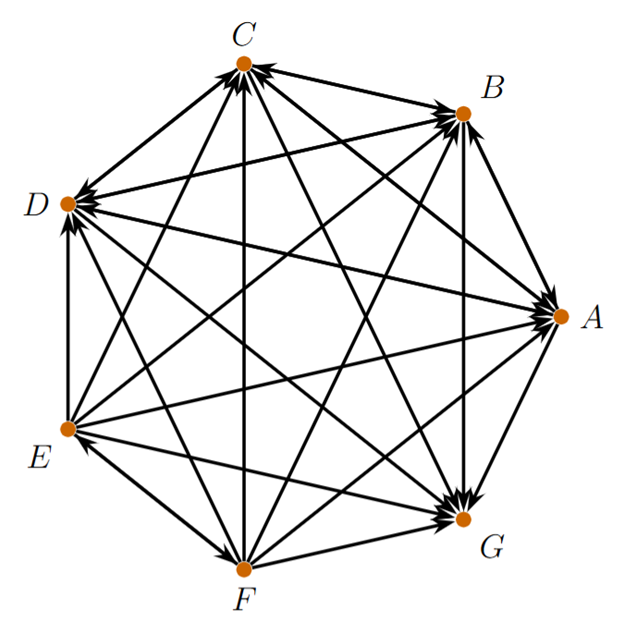}
\caption{Heptagon with microscale directed connectivity graph between the seven blocks $A,B,C,D,E,F,G$ of the sequence space $Seq(\mathbb{R})$.\label{fig3}}
\end{figure}

Figure~\ref{fig3} presents microscale  directed connectivity graph between the seven blocks $A,\ldots,G$ of the sequence space $\operatorname{Seq}(\mathbb{R})$, where a dark-orange node represents a block and a bold arrow $X\to Y$ indicates that every sequence in $X$ can be mapped into $Y$ by a suitable Hadamard product with some connector sequence. The diagram highlights the rich bidirectional connectivity among blocks $A$--$F$, the mutual link between $E$ and $F$, and the fact that $G$ is a global attractor with no outgoing connections.\par

\begin{rema}\label{rem3.11}
For the building blocks  $\{A,B,C,D,E,F,G\}$ of the sequence space $Seq(\mathbb{R})$ there are 42 distinctive pair connections at macro-level and 28 distinctive pair connections at micro-level, respectively. 
\end{rema}

\section{Discussion}\label{sec4}
\subsection{Summary \& Contributions}\label{sec4.1}
In this work we introduced a finite partition of the sequence space \(Seq(\mathbb{R})\), the set of all real-valued sequences on \(\mathbb{R}\), indexed by the pair of values \((\liminf a_n,\limsup a_n)\), and we supplied constructive closed-form examples representing each of the seven resulting blocks each with the size of continuum. The customary Banach subspaces of the sequence space $Seq(\mathbb{R})$, i.e., $c_{00} \subset \ell^{1} \subset \ell^{2} \subset \cdots \subset \ell^{p} \subset c_{0} \subset c \subset \ell^{\infty}, (1\leq p<+\infty)$ lie entirely within just two of the seven blocks in our partition, namely the blocks $B$ and $G$. Thus, the classical Banach-space viewpoint effectively explores only $2/7 \approx 28.6\%$ of the total number of regions of $Seq(\mathbb{R})$. By contrast, the present work provides constructive information about the remaining five blocks $A,C,D,E,F$, thereby complementing the traditional 
picture with a more comprehensive description of the asymptotic landscape of real-valued sequences. At the macroscale we described an ideal pattern of admissible relations between these blocks, while at the microscale we investigated the connections actually realized via pointwise products of sequences, showing that this finer structure covers about two-thirds of the macroscale possibilities. The central notion of \emph{connection} between sequences and between blocks led to a distinctive configuration with 28 realized relations among the seven blocks. Overall, the analysis singles out the block \(G\) of convergent sequences as a particularly prominent component of the partition, since first, it is the only subspace of the $Seq(\mathbb{R})$ , and, at the microscale it is directly connected-as attractor- to every other block corresponding to infinite (non-convergent) behavior.\par 

\subsection{Comparison of Macroscale Matrix $U$ versus Microscale Matrix $V$.}\label{sec4.2}

We compare the macroscale adjacency matrix \(U\)~(\ref{matrix-u}), which encodes all idealized block--to--block connections, with the microscale matrix \(V\)~(\ref{matrix-v}), which records the connections actually realized by our construction. By design, \(U\) has ones at every off–diagonal position (42 ones and 7 diagonal zeros), while \(V\) has 28 ones. Every one in \(V\) appears at a position where \(U\) also has a one, so in terms of edge sets of directed graphs we have:
\begin{eqnarray}
E(V) \subseteq E(U),    
\end{eqnarray}
 
and thus \(V\) is a subgraph of \(U\); nothing in \(V\) contradicts \(U\), it simply omits some edges.

To quantify this compatibility, we use the standard contingency counts
\[
\begin{aligned}
n_{11} &= \#\{(i,j) : U_{ij}=1,\; V_{ij}=1\} = 28,\\
n_{10} &= \#\{(i,j) : U_{ij}=1,\; V_{ij}=0\} = 14,\\
n_{01} &= \#\{(i,j) : U_{ij}=0,\; V_{ij}=1\} = 0,\\
n_{00} &= \#\{(i,j) : U_{ij}=0,\; V_{ij}=0\} = 7,
\end{aligned}
\]
where the last line corresponds to the diagonal zeros.

First, the “coverage of \(U\) by \(V\)” (recall) is given by:
\begin{eqnarray}
\mathrm{Coverage}(U,V)
  = \frac{n_{11}}{n_{11}+n_{10}}
  = \frac{28}{28+14}
  = \frac{28}{42}
  = \frac{2}{3}
  \approx 66.7\%.    
\end{eqnarray}
 
This is exactly the fraction of potential connections allowed by \(U\) that are actually realized in \(V\).

Second, the “consistency of \(V\) with \(U\)” (precision) is:
\begin{eqnarray}
\mathrm{Consistency}(U,V)
  = \frac{n_{11}}{n_{11}+n_{01}}
  = \frac{28}{28+0}
  = 100\%.    
\end{eqnarray}
 
Thus every edge that appears in \(V\) is permitted by \(U\); there are no forbidden edges.

Third, the Jaccard similarity  of the edge sets is \cite{Levy2025Jaccard,Shibata2012}:
\begin{eqnarray}
\mathrm{JacSim}(U,V)
  = \frac{n_{11}}{n_{11}+n_{10}+n_{01}}
  = \frac{28}{28+14+0}
  = \frac{28}{42}
  = \frac{2}{3}
  \approx 66.7\%,    
\end{eqnarray}
 
which in this particular setting coincides with the coverage because \(E(V) \subseteq E(U)\).

Finally, considering all \(7\times 7 = 49\) entries, the Hamming similarity is \cite{Leskovec2020Hamming,Jamil2023}:

\begin{eqnarray}
\mathrm{HamSim}(U,V)
  = \frac{n_{11}+n_{00}}{49}
  = \frac{28+7}{49}
  = \frac{35}{49}
  \approx 71.4\%.    
\end{eqnarray}
 
Taken together, these measures show that \(V\) provides a dense, though not exhaustive, microscale realization of the macroscale connectivity encoded by \(U\): it respects all global constraints while omitting roughly one third of the possible connections.

\subsection{Limitations \& Future Work}\label{sec4.3}
The limitations of this work are transparent and, at the same time, suggest several directions for further investigation. First, our partition of the sequence space \(Seq(\mathbb{R})\) relies solely on the limit profile manifested by the pair \((\liminf a_n,\limsup a_n)\), leaving aside other potentially informative features such as rates of convergence or divergence, oscillatory behavior, periodicity, or the Cauchy property. Second, the connection relation introduced in Definition~\ref{Def3.7}, which underlies the microscale relationships between blocks, is not an equivalence relation and is therefore only a partial tool for organizing these interactions. Third, in the present work we have restricted our attention to the Hadamard (pointwise) product of sequences. There are,  however, several other natural products on \(Seq(\mathbb{R})\), such as the \emph{Cauchy product}, the \emph{Dirichlet product}, and the more 
general \emph{discrete convolution} of sequences. It would be an interesting direction for future research to revisit the  constructions and results developed here under these alternative products, and to investigate how the algebraic and 
topological properties of the sequence blocks change when Cauchy, Dirichlet, or convolution products are used in place of the Hadamard product. Finally, it would be natural to study how the asymptotic equivalence relation of Definition~\ref{Def3.1} shapes both the structure and the number of blocks in our partition, and how this picture changes when the defining ingredients in Definition~\ref{Def3.1} are replaced by alternative properties of real-valued sequences. Such modifications are likely to produce new partitions of \(Seq(\mathbb{R})\) with different block structures, opening a range of promising avenues for a more refined representation of the sequence space.\par

\subsection{Conclusion}\label{sec4.4}
In conclusion, this work provides a constructive, finite, partition-based description of the sequence space \(Seq(\mathbb{R})\), which is more commonly viewed through the lens of its infinite-dimensional linear structure and associated bases. By organizing real-valued sequences into a small number of asymptotic blocks and supplying explicit representatives for each class, we complement the traditional vector-space perspective with a coarse, yet informative, qualitative classification. In addition, the notion of \emph{connection} between sequences and between blocks enriches this picture by encoding how elements of different classes can interact via pointwise products, thereby endowing \(Seq(\mathbb{R})\) with an additional layer of structural insight beyond its standard linear-algebraic description.\par

\subsection*{Funding} This research received no external funding.
\subsection*{Institutional Review Board Statement} Not Applicable.
\subsection*{Informed Consent Statement} Not Applicable.
\subsection*{Data Availability Statement} Not Applicable.
\subsection*{Acknowledgments} Not Applicable.
\subsection*{Conflicts of Interest} The author  declares no conflicts of interest.
\subsection*{Abbreviations}
The following abbreviations are used in this manuscript:\newline
c: cardinality of the continuum; HamSim: Hamming similarity; JacSim: Jaccard similarity; lim inf: limit inferior; lim sup: limit superior; N: set of natural numbers; R: set of real numbers; Seq(R): sequence space.

\end{document}